\documentclass[reqno, 11pt]{amsart}
\usepackage{amssymb}
\usepackage{amsfonts}
\usepackage{srcltx}
\usepackage{enumerate}
\usepackage{cite}

\newtheorem{theorem}{Theorem}[section]
\newtheorem{lemma}[theorem]{Lemma}

\newtheorem{corollary}[theorem]{Corollary}
\theoremstyle{definition}

\theoremstyle{remark}

\DeclareMathOperator{\diam}{diam}
\DeclareMathOperator{\conv}{conv}

\DeclareMathOperator{\Fix}{Fix}
\DeclareMathOperator{\WCS}{WCS}
\DeclareMathOperator{\Nor}{N}

\begin{document}
\title[Asymptotically regular semigroups]{On the structure of fixed-point
sets of asymptotically regular semigroups}
\author[A. Wi\'{s}nicki]{Andrzej Wi\'{s}nicki}

\begin{abstract}
We extend a few recent results of G\'{o}rnicki (2011) asserting that the set
of fixed points of an asymptotically regular mapping is a retract of its
domain. In particular, we prove that in some cases the resulting retraction is
H\"{o}lder continuous. We also characterise Bynum's coefficients and the
Opial modulus in terms of nets.
\end{abstract}

\address{Andrzej Wi\'{s}nicki, Institute of Mathematics, Maria Curie-Sk\l %
odowska University, 20-031 Lublin, Poland}
\email{awisnic@hektor.umcs.lublin.pl}
\maketitle

\section{Introduction.}

The notion of asymptotic regularity, introduced by Browder and Petryshyn in
\cite{BrPe}, has become a standing assumption in many results concerning
fixed points of nonexpansive and more general mappings. Recall that a
mapping $T:M\rightarrow M$ acting on a metric space $(M,d)$ is said to be
asymptotically regular if
\begin{equation*}
\lim_{n\rightarrow \infty }d(T^{n}x,T^{n+1}x)=0
\end{equation*}%
for all $x\in M.$ Ishikawa\ \cite{Is} proved that if $C$ is a bounded closed
convex subset of a Banach space $X$ and $T:C\rightarrow C$ is nonexpansive,
then the mapping $T_{\lambda }=(1-\lambda )I+\lambda T$ is asymptotically
regular for each $\lambda \in (0,1).$ Edelstein and O'Brien \cite{EdOb}
showed independently that $T_{\lambda }$ is uniformly asymptotically regular
over $x\in C,$ and Goebel and Kirk \cite{GoKi3} proved that the convergence
is even uniform with respect to all nonexpansive mappings from $C$ into $C$.
Other examples of asymptotically regular mappings are given by the result of
Anzai and Ishikawa \cite{AnIs} (see also \cite{XuYa}): if $T$ is an affine
mapping acting on a bounded closed convex subset of a locally convex space $%
X $, then $T_{\lambda }=(1-\lambda )I+\lambda T$ is uniformly asymptotically
regular.

In 1987, Lin \cite{Li} constructed a uniformly asymptotically regular
Lipschitz mapping in $\ell _{2}$ without fixed points which extended an
earlier construction of Tingley \cite{Ti}. Subsequently, Maluta, Prus and Wo%
\'{s}ko \cite{MaPrWo} proved that there exists a continuous fixed-point free
asymptotically regular mapping defined on any bounded convex subset of a
normed space which is not totally bounded (see also \cite{Er}). For the
fixed-point existence theorems for asymptotically regular mappings we refer
the reader to the papers by T. Dom\'{\i}nguez Benavides, J. G\'{o}rnicki,
M. A. Jap\'{o}n Pineda and H. K. Xu (see \cite{DoJa, DoXu, Go1}).

It was shown in \cite{SeWi} that the set of fixed points of a k-uniformly
Lipschitzian mapping in a uniformly convex space is a retract of its domain
if $k$ is close to $1$. In recent papers \cite{GoT,GoTai,GoN}, J. G\'{o}rnicki
proved several results concerning the structure of fixed-point sets
of asymptotically regular mappings in uniformly convex spaces. In this paper
we continue this work and extend a few of G\'{o}rnicki's results in two
aspects: we consider a more general class of spaces and prove that in some
cases, the fixed-point set $\mathrm{Fix\,}T$ is not only a (continuous)
retract but even a H\"{o}lder continuous retract of the domain. We present
our results in a more general case of a one-parameter nonlinear semigroup.
We also characterise Bynum's coefficients and the Opial modulus in terms of
nets.

\section{Preliminaries}

Let $G$ be an unbounded subset of $[0,\infty )$ such that $t+s,t-s\in G$ for
all $t,s\in G$ with $t>s$ (e.g., $G=[0,\infty )$ or $G=\mathbb{N}$). By a
nonlinear semigroup on $C$ we shall mean a one-parameter family of mappings $%
\mathcal{T}=\{T_{t}:t\in G\}$ from $C$ into $C$ such that $%
T_{t+s}x=T_{t}\,T_{s}x$ for all $t,s\in G$ and $x\in C$. In particular, we
do not assume in this paper that $\{T_{t}:t\in G\}$ is strongly continuous.
We use a symbol $|T|$ to denote the exact Lipschitz constant of a mapping $%
T:C\rightarrow C$, i.e.,
\begin{equation*}
|T|=\inf \{k:\Vert Tx-Ty\Vert \leq k\Vert x-y\Vert \ \text{for\ all}\ x,y\in
C\}.
\end{equation*}%
If $T$ is not Lipschitzian we define $|T|=\infty $.

A semigroup $\mathcal{T}=\{T_{t}:t\in G\}$ from $C$ into $C$ is said to be
asymptotically regular if $\lim_{t}\left\Vert T_{t+h}x-T_{t}x\right\Vert =0$
for every $x\in C$ and $h\in G.$

Assume now that $C$ is convex and weakly compact and $\mathcal{T}%
=\{T_{t}:t\in G\}$ is a nonlinear semigroup on $C$ such that $s(\mathcal{T}%
)=\liminf_{t}|T_{t}|<\infty .$ Choose a sequence $(t_{n})$ of elements in $G$
such that $\lim_{n\rightarrow \infty }t_{n}=\infty $ and $s(\mathcal{T}%
)=\lim_{n\rightarrow \infty }\left\vert T_{t_{n}}\right\vert .$ By
Tikhonov's theorem, there exists a pointwise weakly convergent subnet $%
(T_{t_{n_{\alpha }}})_{\alpha \in \emph{A}}$ of $(T_{t_{n}}).$ We denote it
briefly by $(T_{t_{\alpha }})_{\alpha \in \emph{A}}.$ For every $x\in C$,
define
\begin{equation}
Lx=w\text{-}\lim_{\alpha }T_{t_{\alpha }}x,  \label{Lx}
\end{equation}%
i.e., $Lx$ is the weak limit of the net $(T_{t_{\alpha }}x)_{\alpha \in
\emph{A}}$. Notice that $Lx$ belongs to $C$ since $C$ is convex and weakly
compact. The weak lower semicontinuity of the norm implies
\begin{equation*}
\Vert Lx-Ly\Vert \leq \liminf_{\alpha }\Vert T_{t_{\alpha }}x-T_{t_{\alpha
}}y\Vert \leq \limsup_{n\rightarrow \infty }\Vert T_{t_{n}}x-T_{t_{n}}y\Vert
\leq s(\mathcal{T})\Vert x-y\Vert .
\end{equation*}%
We formulate the above observation as a separate lemma.

\begin{lemma}
\label{nonexp}Let $C$ be a convex weakly compact subset of a Banach space $X$
and let $\mathcal{T}=\{T_{t}:t\in G\}$ be a semigroup on $C$ such that $s(%
\mathcal{T})=\liminf_{t}|T_{t}|<\infty .$ Then the mapping $L:C\rightarrow C$
defined by (\ref{Lx}) is $s(\mathcal{T})$-Lipschitz.
\end{lemma}

We end this section with the following variant of a well known result which
is crucial for our work (see, e.g., \cite[Prop. 1.10]{BeLi}).

\begin{lemma}
\label{holder}Let $(X,d)$ be a complete bounded metric space and let $%
L:X\rightarrow X$ be a k-Lipschitz mapping. Suppose there exists $0<\gamma
<1 $ and $c>0$ such that $\Vert L^{n+1}x-L^{n}x\Vert \leq c\gamma ^{n}$ for
every $x\in X$. Then $Rx=\lim_{n\rightarrow \infty }L^{n}x$ is a H\"{o}lder
continuous mapping.
\end{lemma}

\begin{proof}
We may assume that $\diam X<1$. Fix $x\neq y$ in $X$ and notice that for any
$n\in \mathbb{N}$,
\begin{equation*}
d(Rx,Ry)\leq d(Rx,L^{n}x)+d(L^{n}x,L^{n}y)+d(L^{n}y,Ry)\leq 2c\frac{\gamma
^{n}}{1-\gamma }+k^{n}d(x,y).
\end{equation*}%
Take $\alpha <1$ such that $k\leq \gamma ^{1-\alpha ^{-1}}$ and put $\gamma
^{n-r}=d(x,y)^{\alpha }$ for some $n\in \mathbb{N}$ and $0<r\leq 1$. Then $%
k^{n-1}\leq (\gamma ^{1-\alpha ^{-1}})^{n-r}$ and hence
\begin{equation*}
d(Rx,Ry)\leq 2c\frac{\gamma ^{n-r}}{1-\gamma }+k(\gamma ^{n-r})^{1-\alpha
^{-1}}d(x,y)=(\frac{2c}{1-\gamma }+k)d(x,y)^{\alpha }.
\end{equation*}
\end{proof}

\section{Bynum's coefficients and Opial's modulus in terms of nets}

From now on, $C$ denotes a nonempty convex weakly compact subset of a Banach
space $X$. Let $\mathcal{A}$ be a directed set, $(x_{\alpha })_{\alpha \in
\mathcal{A}}$ a bounded net in $X$, $y\in X$ and write
\begin{align*}
r(y,(x_{\alpha }))& =\limsup_{\alpha }\Vert x_{\alpha }-y\Vert , \\
r(C,(x_{\alpha }))& =\inf \{r(y,(x_{\alpha })):y\in C\}, \\
A(C,(x_{\alpha }))& =\{y\in C:r(y,(x_{\alpha }))=r(C,(x_{\alpha }))\}.
\end{align*}

The number $r(C,(x_{\alpha }))$ and the set $A(C,(x_{\alpha }))$ are called,
respectively, the asymptotic radius and the asymptotic center of $(x_{\alpha
})_{\alpha \in \mathcal{A}}$ relative to $C$. Notice that $A(C,(x_{\alpha
})) $ is nonempty convex and weakly compact. Write
\begin{equation*}
r_{a}(x_{\alpha })=\inf \{\limsup_{\alpha }\Vert x_{\alpha }-y\Vert :y\in
\overline{\conv}(\{x_{\alpha }:\alpha \in \mathcal{A}\})\}
\end{equation*}%
and let
\begin{equation*}
\diam_{a}(x_{\alpha })=\inf_{\alpha }\sup_{\beta ,\gamma \geq \alpha }\Vert
x_{\beta }-x_{\gamma }\Vert
\end{equation*}%
denote the asymptotic diameter of $(x_{\alpha })$.

The normal structure coefficient $\Nor(X)$ of a Banach space $X$ is defined
by
\begin{equation*}
\Nor(X)=\sup \left\{ k:k\,r(K)\leq \diam K\ \ \text{for\ each\ bounded\
convex\ set}\ K\subset X\right\} ,
\end{equation*}%
where $r(K)=\inf_{y\in K}\sup_{x\in K}\Vert x-y\Vert $ is the Chebyshev
radius of $K$ relative to itself. Assuming that $X$ does not have the Schur
property, the weakly convergent sequence coefficient (or Bynum's
coefficient) is given by
\begin{equation*}
\WCS(X)%
=\sup \left\{ k:k\,r_{a}(x_{n})\leq \diam_{a}(x_{n})\ \ \text{for\ each\
sequence}\ x_{n}\overset{w}{\longrightarrow }0\right\} ,
\end{equation*}%
where $x_{n}\overset{w}{\longrightarrow }0$ means that $(x_{n})$ is weakly
null in $X$ (see \cite{By}). For Schur spaces, we define $WCS(X)=2$.

It was proved independently in \cite{DoLop, Pr, Zh} that

\begin{equation}
WCS(X)=\sup \left\{ k:k\,\limsup_{n}\Vert x_{n}\Vert \leq \diam_{a}(x_{n})\
\text{for each\ sequence}\ x_{n}\overset{w}{\longrightarrow }0\right\}
\label{wcs1}
\end{equation}%
and, in \cite{DoLoXu}, that
\begin{equation*}
WCS(X)=\sup \left\{ k:k\,\limsup_{n}\Vert x_{n}\Vert \leq D[(x_{n})]\ \text{%
for\ each\ sequence}\ x_{n}\overset{w}{\longrightarrow }0\right\} ,
\end{equation*}%
where $D[(x_{n})]=\limsup_{m}\limsup_{n}\left\Vert x_{n}-x_{m}\right\Vert
.\smallskip $

Kaczor and Prus \cite{KaPr} initiated a systematic study of assumptions
under which one can replace sequences by nets in a given condition. We
follow the arguments from that paper and use the well known method of
constructing basic sequences attributed to S. Mazur (see \cite{Pe}). Let us
first recall a variant of a classical lemma which can be proved in the same
way as for sequences (see, e.g., \cite[Lemma]{Pe}).

\begin{lemma}
\label{Ma} Let $\{x_{\alpha }\}_{\alpha \in \mathcal{A}}$ be a bounded net
in $X$ weakly converging to $0$ such that $\inf_{\alpha }\Vert x_{\alpha
}\Vert >0$. Then for every $\varepsilon >0$, $\alpha ^{\prime }\in \mathcal{A%
}$ and for every finite dimensional subspace $E$ of $X$, there is $\alpha
>\alpha ^{\prime }$ such that
\begin{equation*}
\Vert e+tx_{\alpha }\Vert \geq (1-\varepsilon )\Vert e\Vert
\end{equation*}%
for any $e\in E$ and every scalar $t.$
\end{lemma}

Recall that a sequence $(x_{n})$ is basic if and only if there exists a
number $c>0$ such that $\Vert \sum_{i=1}^{q}t_{i}x_{i}\Vert \leq c\Vert
\sum_{i=1}^{p}t_{i}x_{i}\Vert $ for any integers $p>q\geq 1$ and any
sequence of scalars $(t_{i})$. In the proof of the next lemma, based on
Mazur's technique, we follow in part the reasoning given in \cite[Cor. 2.6]%
{KaPr}. Set $D[(x_{\alpha })]=\limsup_{\alpha }\limsup_{\beta }\left\Vert
x_{\alpha }-x_{\beta }\right\Vert .$

\begin{lemma}
\label{KaPr} Let $(x_{\alpha })_{\alpha \in \mathcal{A}}$ be a bounded net
in $X$ which converges to $0$ weakly but not in norm. Then there exists an
increasing sequence $(\alpha _{n})$ of elements of $\mathcal{A}$ such that $%
\lim_{n}\Vert x_{\alpha _{n}}\Vert =\limsup_{\alpha }\Vert x_{\alpha }\Vert $%
, $\diam_{a}(x_{\alpha _{n}})\leq D[(x_{\alpha })]$ and $(x_{\alpha _{n}})$
is a basic sequence.
\end{lemma}

\begin{proof}
Since $(x_{\alpha })_{\alpha \in \mathcal{A}}$ does not converge strongly to
$0$ and $D[(x_{\alpha _{s}})]\leq D[(x_{\alpha })]$ for any subnet $%
(x_{\alpha _{s}})_{s\in \mathcal{B}}$ of $(x_{\alpha })_{\alpha \in \mathcal{%
A}}$, we can assume, passing to a subnet, that $\inf_{\alpha }\Vert
x_{\alpha }\Vert >0$ and the limit $c=\lim_{\alpha }\Vert x_{\alpha }\Vert $
exists. Write $d=D[(x_{\alpha })]$. Let $(\varepsilon _{n})$ be a sequence
of reals from the interval $(0,1)$ such that $\Pi _{n=1}^{\infty
}(1-\varepsilon _{n})>0$. We shall define the following sequences $(\alpha
_{n})$ and $(\beta _{n})$ by induction.

Let us put $\alpha _{1}<\beta _{1}\in \mathcal{A}$ such that $\left\vert
\Vert x_{\alpha _{1}}\Vert -c\right\vert <1$ and $\sup_{\beta \geq \beta
_{1}}\Vert x_{\alpha _{1}}-x_{\beta }\Vert <d+1$. By the definitions of $c$
and $d$, there exists $\alpha ^{\prime }>\beta _{1}$ such that $\left\vert
\Vert x_{\alpha }\Vert -c\right\vert <\frac{1}{2}$ and $\inf_{\beta ^{\prime
}}\sup_{\beta \geq \beta ^{\prime }}\Vert x_{\alpha }-x_{\beta }\Vert <d+%
\frac{1}{2}$ for every $\alpha \geq \alpha ^{\prime }.$ It follows from
Lemma \ref{Ma} that there exists $\alpha _{2}>\alpha ^{\prime }$ such that%
\begin{equation*}
\Vert t_{1}x_{\alpha _{1}}+t_{2}x_{\alpha _{2}}\Vert \geq (1-\varepsilon
_{2})\Vert t_{1}x_{\alpha _{1}}\Vert
\end{equation*}%
for any scalars $t_{1},t_{2}.$ Furthermore, $\left\vert \Vert x_{\alpha
_{2}}\Vert -c\right\vert <\frac{1}{2},$ and we can find $\beta _{2}>\alpha
_{2}$ such that $\sup_{\beta \geq \beta _{2}}\Vert x_{\alpha _{2}}-x_{\beta
}\Vert <d+\frac{1}{2}.$

Suppose now that we have chosen $\alpha _{1}<\beta _{1}<...<\alpha
_{n}<\beta _{n}$ $(n>1)$ in such a way that $\left\vert \Vert x_{\alpha
_{k}}\Vert -c\right\vert <\frac{1}{k}$, $\sup_{\beta \geq \beta _{k}}\Vert
x_{\alpha _{k}}-x_{\beta }\Vert <d+\frac{1}{k}$ and
\begin{equation*}
(1-\varepsilon _{k})\Vert t_{1}x_{\alpha _{1}}+...+t_{k-1}x_{\alpha
_{k-1}}\Vert \leq \Vert t_{1}x_{\alpha _{1}}+...+t_{k}x_{\alpha _{k}}\Vert
\end{equation*}%
for any scalars $t_{1},...,t_{k}$, $k=2,...,n.$ From the definitions of $c$
and $d$, and by Lemma \ref{Ma}, we can find $\beta _{n+1}>\alpha
_{n+1}>\beta _{n}$ such that $\left\vert \Vert x_{\alpha _{n+1}}\Vert
-c\right\vert <\frac{1}{n+1}$, $\sup_{\beta \geq \beta _{n+1}}\Vert
x_{\alpha _{n+1}}-x_{\beta }\Vert \leq d+\frac{1}{n+1}$ and (considering a
subspace $E$ spanned by the elements $x_{\alpha _{1}},...,x_{\alpha _{n}}$
and putting $e=t_{1}x_{\alpha _{1}}+...+t_{n}x_{\alpha _{n}}$),
\begin{equation*}
(1-\varepsilon _{n+1})\Vert t_{1}x_{\alpha _{1}}+...+t_{n}x_{\alpha
_{n}}\Vert \leq \Vert t_{1}x_{\alpha _{1}}+...+t_{n+1}x_{\alpha _{n+1}}\Vert
\end{equation*}%
for any scalars $t_{1},...,t_{n+1}$.

Notice that the sequence $(x_{\alpha _{n}})$ defined in this way satisfies $%
\lim_{n\rightarrow \infty }\Vert x_{\alpha _{n}}\Vert =c$ and $\diam%
_{a}(x_{\alpha _{n}})\leq d$. Furthermore,
\begin{equation*}
\Vert t_{1}x_{\alpha _{1}}+...+t_{p}x_{\alpha _{p}}\Vert \geq \Pi
_{n=q+1}^{p}(1-\varepsilon _{n})\Vert t_{1}x_{\alpha
_{1}}+...+t_{q}x_{\alpha _{q}}\Vert
\end{equation*}%
for any integers $p>q\geq 1$ and any sequence of scalars $(t_{i})$. Hence $%
(x_{\alpha _{n}})$ is a basic sequence.
\end{proof}

We are now in a position to give a characterization of the coefficient
WCS(X) in terms of nets. The abbreviation \textquotedblleft $\left\{
x_{\alpha }\right\} $ is r.w.c.\textquotedblright\ means that the set $%
\left\{ x_{\alpha }:\alpha \in \mathcal{A}\right\} $ is relatively weakly
compact.

\begin{theorem}
\label{Wi1} Let $X$ be a Banach space without the Schur property and write%
\begin{align*}
w_{1}& =\sup \left\{ k:k\,r_{a}(x_{\alpha })\leq \diam_{a}(x_{\alpha })\
\text{ for\ each\ net}\ x_{\alpha }\overset{w}{\longrightarrow }0,\text{ }%
\left\{ x_{\alpha }\right\} \text{ is r.w.c.}\right\} , \\
w_{2}& =\sup \left\{ k:k\,\limsup_{\alpha }\Vert x_{\alpha }\Vert \leq \diam%
_{a}(x_{\alpha })\ \text{for\ each\ net}\ x_{\alpha }\overset{w}{%
\longrightarrow }0,\text{ }\left\{ x_{\alpha }\right\} \text{ is r.w.c.}%
\right\} , \\
w_{3}& =\sup \left\{ k:k\,\limsup_{\alpha }\Vert x_{\alpha }\Vert \leq
D[(x_{\alpha })]\ \text{for\ each\ net}\ x_{\alpha }\overset{w}{%
\longrightarrow }0,\text{ }\left\{ x_{\alpha }\right\} \text{ is r.w.c.}%
\right\} .
\end{align*}%
Then
\begin{equation*}
\WCS(X)=w_{1}=w_{2}=w_{3}.
\end{equation*}
\end{theorem}

\begin{proof}
Fix $k>w_{3}$ and choose a weakly null net $(x_{\alpha })$ such that the set
$\left\{ x_{\alpha }:\alpha \in \mathcal{A}\right\} $ is relatively weakly
compact and $k\,\limsup_{\alpha }\Vert x_{\alpha }\Vert >D[(x_{\alpha })].$
Then, by Lemma \ref{KaPr}, there exists an increasing sequence $(\alpha
_{n}) $ such that
\begin{equation*}
k\,\lim_{n}\Vert x_{\alpha _{n}}\Vert >D[(x_{\alpha })]\geq \diam%
_{a}(x_{\alpha _{n}})
\end{equation*}
and $(x_{\alpha _{n}})$ is a basic sequence. Since the set $\left\{
x_{\alpha }:\alpha \in \mathcal{A}\right\} $ is relatively weakly compact,
we can assume (passing to a subsequence) that $(x_{\alpha _{n}})$ is weakly
convergent. Since it is a basic sequence, its weak limit equals zero. It
follows from (\ref{wcs1}) that $%
\WCS(X)%
\leq k$ and letting $k$ go to $w_{3}$ we have
\begin{equation*}
\WCS(X)%
\leq w_{3}\leq w_{2}\leq w_{1}\leq
\WCS(X)%
.
\end{equation*}
\end{proof}

Notice that a similar characterisation holds for the normal structure
coefficient.

\begin{theorem}
For a Banach space $X$,%
\begin{equation*}
\Nor(X)=\sup \left\{ k:k\,r_{a}(x_{\alpha })\leq \diam_{a}(x_{\alpha })\
\text{for\ each bounded net}\ (x_{\alpha })\text{ in }X\right\} .
\end{equation*}
\end{theorem}

\begin{proof}
Let
\begin{equation*}
N_{1}=\sup \left\{ k:k\,r_{a}(x_{\alpha })\leq \diam_{a}(x_{\alpha })\ \text{%
for\ each bounded net}\ (x_{\alpha })\text{ in }X\right\} .
\end{equation*}%
Set $k>N_{1}$ and choose a bounded net $(x_{\alpha })$ such that $%
k\,r_{a}(x_{\alpha })>\diam_{a}(x_{\alpha }).$ Fix $y\in \overline{\conv}%
(\{x_{\alpha }:\alpha \in \mathcal{A}\})$ and notice that $%
k\,\limsup_{\alpha }\Vert x_{\alpha }-y\Vert >\diam_{a}(x_{\alpha }).$ In a
straightforward way, we can choose a sequence $(\alpha _{n})$ such that
\begin{equation*}
k\,\lim_{n}\Vert x_{\alpha _{n}}-y\Vert =k\,\limsup_{\alpha }\Vert x_{\alpha
}-y\Vert >\diam_{a}(x_{\alpha })\geq \diam_{a}(x_{\alpha _{n}}).
\end{equation*}%
It follows from \cite[Th. 1]{By} that $\Nor(X)\leq k$ and letting $k$ go to $%
N_{1}$ we have $\Nor(X)\leq N_{1}.$ By \cite[Th. 1]{Lim}, $\Nor(X)\geq N_{1}$
and the proof is complete.
\end{proof}

In the next section we shall need a similar characterisation for the Opial
modulus of a Banach space $X,$ defined for each $c\geq 0$ by
\begin{equation*}
r_{X}(c)=\inf \left\{ \liminf_{n\rightarrow \infty }\left\Vert
x_{n}+x\right\Vert -1\right\} ,
\end{equation*}%
where the infimum is taken over all $x\in X$ with $\left\Vert x\right\Vert
\geq c$ and all weakly null sequences $(x_{n})$ in $X$ such that $%
\liminf_{n\rightarrow \infty }\left\Vert x_{n}\right\Vert \geq 1$ (see \cite%
{LiTaXu}). We first prove the following counterpart of Lemma \ref{KaPr}.

\begin{lemma}
\label{KaPr2} Let $(x_{\alpha })_{\alpha \in \mathcal{A}}$ be a bounded net
in $X$ which converges to $0$ weakly but not in norm and $x\in X.$ Then
there exists an increasing sequence $(\alpha _{n})$ of elements of $\mathcal{%
A}$ such that $\lim_{n}\Vert x_{\alpha _{n}}+x\Vert =\liminf_{\alpha }\Vert
x_{\alpha }+x\Vert ,$ $\lim_{n}\Vert x_{\alpha _{n}}\Vert \geq
\liminf_{\alpha }\Vert x_{\alpha }\Vert $ and $(x_{\alpha _{n}})$ is a basic
sequence.
\end{lemma}

\begin{proof}
Since $(x_{\alpha })_{\alpha \in \mathcal{A}}$ does not converge strongly to
$0$ and
\begin{equation*}
\liminf_{s}\Vert x_{\alpha _{s}}\Vert \geq \liminf_{\alpha }\Vert x_{\alpha
}\Vert
\end{equation*}%
for any subnet $(x_{\alpha _{s}})_{s\in \mathcal{B}}$ of $(x_{\alpha
})_{\alpha \in \mathcal{A}}$, it is sufficient (passing to a subnet) to
consider only the case that $\inf_{\alpha }\Vert x_{\alpha }\Vert >0$ and
the limits $c_{1}=\liminf_{\alpha }\Vert x_{\alpha }+x\Vert $, $%
c_{2}=\liminf_{\alpha }\Vert x_{\alpha }\Vert $ exist. Let $(\varepsilon
_{n})$ be a sequence of reals from the interval $(0,1)$ such that $\Pi
_{n=1}^{\infty }(1-\varepsilon _{n})>0$. We shall define the sequence $%
(\alpha _{n})$ by induction.

Let us put $\alpha _{1}\in \mathcal{A}$ such that $\left\vert \Vert
x_{\alpha _{1}}+x\Vert -c_{1}\right\vert <1$ and $\left\vert \Vert x_{\alpha
_{1}}\Vert -c_{2}\right\vert <1$. By the definitions of $c_{1}$ and $c_{2}$,
there exists $\alpha ^{\prime }>\alpha _{1}$ such that $\left\vert \Vert
x_{\alpha }+x\Vert -c_{1}\right\vert <\frac{1}{2}$ and $\left\vert \Vert
x_{\alpha }\Vert -c_{2}\right\vert <\frac{1}{2}$ for every $\alpha \geq
\alpha ^{\prime }.$ It follows from Lemma \ref{Ma} that there exists $\alpha
_{2}>\alpha ^{\prime }$ such that%
\begin{equation*}
\Vert t_{1}x_{\alpha _{1}}+t_{2}x_{\alpha _{2}}\Vert \geq (1-\varepsilon
_{2})\Vert t_{1}x_{\alpha _{1}}\Vert
\end{equation*}%
for any scalars $t_{1},t_{2}.$ We can now proceed analogously to the proof
of Lemma \ref{KaPr} to obtain a basic sequence $(x_{\alpha _{n}})$ with the
desired properties.
\end{proof}

\begin{theorem}
\label{Wi2}For a Banach space $X$ without the Schur property and for $c\geq
0,$%
\begin{equation*}
r_{X}(c)=\inf \left\{ \liminf_{\alpha }\left\Vert x_{\alpha }+x\right\Vert
-1\right\} ,
\end{equation*}%
where the infimum is taken over all $x\in X$ with $\left\Vert x\right\Vert
\geq c$ and all weakly null nets $(x_{\alpha })$ in $X$ such that $%
\liminf_{\alpha }\Vert x_{\alpha }\Vert \geq 1$ and the set $\left\{
x_{\alpha }:\alpha \in \mathcal{A}\right\} $ is relatively weakly compact.
\end{theorem}

\begin{proof}
Let $r_{1}(c)=\inf \left\{ \liminf_{\alpha }\left\Vert x_{\alpha
}+x\right\Vert -1\right\} ,$ where the infimum is taken as above. Fix $c\geq
0$ and take $k>r_{1}(c).$ Then there exist $x\in X$ with $\left\Vert
x\right\Vert \geq c$ and a weakly null net $(x_{\alpha })_{\alpha \in
\mathcal{A}}$ such that $\liminf_{\alpha }\Vert x_{\alpha }\Vert \geq 1,$ $%
\left\{ x_{\alpha }:\alpha \in \mathcal{A}\right\} $ is relatively weakly
compact and
\begin{equation*}
\liminf_{\alpha }\left\Vert x_{\alpha }+x\right\Vert -1<k.
\end{equation*}%
By Lemma \ref{KaPr2}, there exists an increasing sequence $(\alpha _{n})$ of
elements of $\mathcal{A}$ such that $\lim_{n}\Vert x_{\alpha _{n}}\Vert \geq
1,\lim_{n}\Vert x_{\alpha _{n}}+x\Vert -1<k$ and $(x_{\alpha _{n}})$ is a
basic sequence. Since $\left\{ x_{\alpha }:\alpha \in \mathcal{A}\right\} $
is relatively weakly compact, we can assume (passing to a subsequence) that $%
(x_{\alpha _{n}})$ is weakly null. Hence $r_{X}(c)<k$ and since $k$ is an
arbitrary number greater than $r_{1}(c)$, it follows that $r_{X}(c)\leq
r_{1}(c).$ The reverse inequality is obvious.
\end{proof}

\section{Fixed-point sets as H\"{o}lder continuous retracts}

The following lemma may be proved in a similar way to \cite[Th. 7.2 ]{DoJaLo}%
.

\begin{lemma}
\label{main}Let $C$ be a nonempty convex weakly compact subset of a Banach
space $X$ and $\mathcal{T}=\{T_{t}:t\in G\}$ an asymptotically regular
semigroup on $C$ such that $s(\mathcal{T})=\lim_{\alpha }\left\vert
T_{t_{\alpha }}\right\vert $ for a pointwise weakly convergent subnet $%
(T_{t_{\alpha }})_{\alpha \in \emph{A}}$ of $(T_{t})_{t\in G}.$ Let $%
x_{0}\in C$, $x_{m+1}=w$-$\lim_{\alpha }T_{t_{\alpha }}x_{m},m=0,1,...,$ and
\begin{equation*}
B_{m}=\limsup_{\alpha }\left\Vert T_{t_{\alpha }}x_{m}-x_{m+1}\right\Vert .
\end{equation*}%
Assume that
\end{lemma}

\begin{enumerate}
\item[(a)] $s(\mathcal{T})<\sqrt{\WCS(X)}$ or,

\smallskip

\item[(b)] $s(\mathcal{T})<1+r_{X}(1).$
\end{enumerate}
\smallskip Then, there exists $\gamma <1$ such that $B_{m}\leq \gamma
B_{m-1} $ for any $m=1,2,...$.

\begin{proof}
It follows from the asymptotic regularity of $\{T_{t}:t\in G\}$ that
\begin{equation*}
\limsup_{\alpha }\left\Vert T_{t_{\alpha }-l}\,x-y\right\Vert
=\limsup_{\alpha }\left\Vert T_{t_{\alpha }}x-y\right\Vert
\end{equation*}%
for any $l\in G$ and $x,y\in C$. Thus%
\begin{align*}
& D[(T_{t_{\alpha }}x_{m})]=\limsup_{\beta }\limsup_{\alpha }\left\Vert
T_{t_{\alpha }}x_{m}-T_{t_{\beta }}x_{m}\right\Vert \\
& \ \leq \limsup_{\beta }\left\vert T_{t_{\beta }}\right\vert
\limsup_{\alpha }\left\Vert T_{t_{\alpha }-t_{\beta }}x_{m}-x_{m}\right\Vert
=s(\mathcal{T})\limsup_{\alpha }\left\Vert T_{t_{\alpha
}}x_{m}-x_{m}\right\Vert .
\end{align*}%
Hence, from Theorem \ref{Wi1} and from the weak lower semicontinuity of the
norm,%
\begin{align*}
B_{m}& \leq \frac{D[(T_{t_{\alpha }}x_{m})]}{\WCS(X)}\leq \frac{s(\mathcal{T}%
)}{\WCS(X)}\limsup_{\alpha }\left\Vert T_{t_{\alpha }}x_{m}-x_{m}\right\Vert
\\
& \leq \frac{s(\mathcal{T})}{\WCS(X)}\limsup_{\alpha }\liminf_{\beta
}\left\Vert T_{t_{\alpha }}x_{m}-T_{t_{\beta }}x_{m-1}\right\Vert \\
& \leq \frac{s(\mathcal{T})}{\WCS(X)}\limsup_{\alpha }\left\vert
T_{t_{\alpha }}\right\vert \limsup_{\beta }\left\Vert x_{m}-T_{t_{\beta
}-t_{\alpha }}x_{m-1}\right\Vert =\frac{(s(\mathcal{T}))^{2}}{\WCS(X)}%
B_{m-1}.
\end{align*}%
This gives (a). For (b), we can use Theorem \ref{Wi2} and proceed
analogously to the proof of \cite[Th. 7.2 ]{DoJaLo} (see also \cite[Th. 5]%
{GoN}).
\end{proof}

We are now in a position to prove a qualitative semigroup version of \cite[%
Th. 7.2 (a) (b)]{DoJaLo} which is in turn based on the results given in \cite%
{DoJa, DoXu} (see also \cite{Ku}). It also extends, in a few directions,
\cite[Th. 5]{GoN}.

\begin{theorem}
\label{Thwcs}Let $C$ be a nonempty convex weakly compact subset of a Banach
space $X$ and $\mathcal{T}=\{T_{t}:t\in G\}$ an asymptotically regular
semigroup on $C.$ Assume that
\end{theorem}

\begin{enumerate}
\item[(a)] $s(\mathcal{T})<\sqrt{\WCS(X)}$ or,

\smallskip

\item[(b)] $s(\mathcal{T})<1+r_{X}(1).$
\end{enumerate}
\smallskip Then $\mathcal{T}$ has a fixed point in $C$ and $\Fix\mathcal{T}%
=\{x\in C:T_{t}x=x,\,t\in G\}$ is a H\"{o}lder continuous retract of $C.$

\begin{proof}
Choose a sequence $(t_{n})$ of elements in $G$ such that $\lim_{n\rightarrow
\infty }t_{n}=\infty $ and $s(\mathcal{T})=\lim_{n\rightarrow \infty
}\left\vert T_{t_{n}}\right\vert .$ Let $(T_{t_{n_{\alpha }}})_{\alpha \in
\emph{A}}$ (denoted briefly by $(T_{t_{\alpha }})_{\alpha \in \emph{A}}$) be
a pointwise weakly convergent subnet of $(T_{t_{n}}).$ Define, for every $%
x\in C$,
\begin{equation*}
Lx=w-\lim_{\alpha }T_{t_{\alpha }}x.
\end{equation*}%
Fix $x_{0}\in C$ and put $x_{m+1}=Lx_{m},m=0,1,....$ Let $%
B_{m}=\limsup_{\alpha }\left\Vert T_{t_{\alpha }}x_{m}-x_{m+1}\right\Vert .$
By Lemma \ref{main}, there exists $\gamma <1$ such that $B_{m}\leq \gamma
B_{m-1}$ for any $m\geq 1.$ Since the norm is weak lower semicontinuous and
the semigroup is asymptotically regular,%
\begin{align*}
& \Vert L^{m+1}x_{0}-L^{m}x_{0}\Vert =\left\Vert x_{m+1}-x_{m}\right\Vert
\leq \liminf_{\alpha }\left\Vert T_{t_{\alpha }}x_{m}-x_{m}\right\Vert \\
& \ \leq \liminf_{\alpha }\liminf_{\beta }\left\Vert T_{t_{\alpha
}}x_{m}-T_{t_{\beta }}x_{m-1}\right\Vert \leq \limsup_{\alpha }\left\vert
T_{t_{\alpha }}\right\vert \limsup_{\beta }\left\Vert x_{m}-T_{t_{\beta
}-t_{\alpha }}x_{m-1}\right\Vert \\
& \ =s(\mathcal{T})B_{m-1}\leq s(\mathcal{T})\gamma ^{m-1}\diam C
\end{align*}%
for every $x_{0}\in C$ and $m\geq 1.$ Furthermore, by Lemma \ref{nonexp},
the mapping $L:C\rightarrow C$ is $s(\mathcal{T})$-Lipschitz. It follows
from Lemma \ref{holder} that $Rx=\lim_{n\rightarrow \infty }L^{n}x$ is a H%
\"{o}lder continuous mapping on $C$. We show that $R$ is a retraction onto $%
\Fix\mathcal{T}.$ It is clear that if $x\in \Fix\mathcal{T},$ then $Rx=x.$
Furthermore, for every $x\in C,m\geq 1$ and $\alpha \in \emph{A},$%
\begin{equation*}
\Vert T_{t_{\alpha }}Rx-Rx\Vert \leq \left\Vert T_{t_{\alpha
}}Rx-T_{t_{\alpha }}L^{m}x\right\Vert +\left\Vert T_{t_{\alpha
}}L^{m}x-L^{m+1}x\right\Vert +\left\Vert L^{m+1}x-Rx\right\Vert
\end{equation*}%
and hence%
\begin{equation*}
\lim_{\alpha }\Vert T_{t_{\alpha }}Rx-Rx\Vert \leq s(\mathcal{T})\left\Vert
Rx-L^{m}x\right\Vert +B_{m}+\left\Vert L^{m+1}x-Rx\right\Vert .
\end{equation*}%
Letting $m$ go to infinity, $\limsup_{\alpha }\Vert T_{t_{\alpha
}}Rx-Rx\Vert =0.$ Since $s(\mathcal{T})=\lim_{\beta }\left\vert T_{t_{\beta
}}\right\vert <\infty ,$ there exists $\beta _{0}\in \emph{A}$ such that $%
\left\vert T_{t_{\beta }}\right\vert <\infty $ for every $\beta \geq \beta
_{0}.$ Then, the asymptotic regularity of $\mathcal{T}$ implies%
\begin{align*}
\Vert T_{t_{\beta }}Rx-Rx\Vert & \leq \left\vert T_{t_{\beta }}\right\vert
\limsup_{\alpha }\Vert Rx-T_{t_{\alpha }}Rx\Vert +\lim_{\alpha }\Vert
T_{t_{\beta }+t_{\alpha }}Rx-T_{t_{\alpha }}Rx\Vert \\
& +\limsup_{\alpha }\Vert T_{t_{\alpha }}Rx-Rx\Vert =0.
\end{align*}%
Hence $T_{t_{\beta }}Rx=Rx$ for every $\beta \geq \beta _{0}$ and, from the
asymptotic regularity again,
\begin{equation*}
\Vert T_{t}Rx-Rx\Vert =\lim_{\beta }\left\Vert T_{t+t_{\beta
}}Rx-T_{t_{\beta }}Rx\right\Vert =0
\end{equation*}%
for each $t\in G.$ Thus $Rx\in \Fix\mathcal{T}$ for every $x\in C$ and the
proof is complete.
\end{proof}

It is well known that the Opial modulus of a Hilbert space $H,$%
\begin{equation*}
r_{H}(c)=\sqrt{1+c^{2}}-1,
\end{equation*}%
and the Opial modulus of $\ell _{p},p>1,$
\begin{equation*}
r_{\ell _{p}}(c)=(1+c^{p})^{1/p}-1
\end{equation*}%
for all $c\geq 0$ (see \cite{LiTaXu}). The following corollaries are
sharpened versions of \cite[Th. 2.2]{GoT} and \cite[Cor. 8]{GoN}.

\begin{corollary}
Let $C$ be a nonempty bounded closed convex subset of a Hilbert space $H.$
If $\mathcal{T}=\{T_{t}:t\in G\}$ is an asymptotically regular semigroup on $%
C$ such that
\begin{equation*}
\liminf_{t}|T_{t}|<\sqrt{2},
\end{equation*}%
then $\Fix\mathcal{T}$ is a H\"{o}lder continuous retract of $C.$
\end{corollary}

\begin{corollary}
Let $C$ be a nonempty bounded closed convex subset of $\ell _{p},1<p<\infty
. $ If $\mathcal{T}=\{T_{t}:t\in G\}$ is an asymptotically regular semigroup
on $C$ such that
\begin{equation*}
\liminf_{t}|T_{t}|<2^{1/p},
\end{equation*}%
then $\Fix\mathcal{T}$ is a H\"{o}lder continuous retract of $C.$
\end{corollary}

Let $1\leq p,q<\infty .$ Recall that the Bynum space $\ell _{p,q}$ is the
space $\ell _{p}$ endowed with the equivalent norm $\Vert x\Vert
_{p,q}=(\Vert x^{+}\Vert _{p}^{q}+\Vert x^{-}\Vert _{p}^{q})^{1/q},$ where $%
x^{+},x^{-}$ denote, respectively, the positive and the negative part of $x.$
If $p>1,$ then
\begin{equation*}
r_{\ell _{p,q}}(c)=\min \{(1+c^{p})^{1/p}-1,(1+c^{q})^{1/q}-1\}
\end{equation*}%
for all $c\geq 0$ (see, e.g., \cite{AyDoLo}). The following corollary
extends \cite[Cor. 10]{GoN}.

\begin{corollary}
Let $C$ be a nonempty convex weakly compact subset of $\ell
_{p,q},1<p<\infty ,1\leq q<\infty .$ If $\mathcal{T}=\{T_{t}:t\in G\}$ is an
asymptotically regular semigroup on $C$ such that
\begin{equation*}
\liminf_{t}|T_{t}|<\min \{2^{1/p},2^{1/q}\},
\end{equation*}%
then $\Fix\mathcal{T}$ is a H\"{o}lder continuous retract of $C.$
\end{corollary}

Let us now examine the case of $p$-uniformly convex spaces. Recall that a
Banach space $X$ is $p$-uniformly convex if $\inf_{\varepsilon >0}\delta
(\varepsilon )\varepsilon ^{-p}>0,$ where $\delta $ denotes the modulus of
uniform convexity of $X.$ If $X$ is $p$-uniformly convex, then (see \cite%
{Xu91})%
\begin{equation}
\left\Vert \lambda x+(1-\lambda )y\right\Vert ^{p}\leq \lambda \left\Vert
x\right\Vert ^{p}+(1-\lambda )\left\Vert y\right\Vert
^{p}-c_{p}W_{p}(\lambda )\left\Vert x-y\right\Vert ^{p}  \label{ineq_Xu}
\end{equation}%
for some $c_{p}>0$ and every $x,y\in X,0\leq \lambda \leq 1,$ where $%
W_{p}(\lambda )=\lambda (1-\lambda )^{p}+\lambda ^{p}(1-\lambda ).$ A Banach
space $X$ satisfies the Opial property if%
\begin{equation*}
\liminf_{n\rightarrow \infty }\left\Vert x_{n}-x\right\Vert
<\liminf_{n\rightarrow \infty }\left\Vert x_{n}-y\right\Vert
\end{equation*}%
for every sequence $x_{n}\overset{w}{\longrightarrow }x$ and $y\neq x.$

The following theorem is an extension of \cite[Th. 7]{GoN}, and a partial
extension of \cite[Th. 9]{GoTai}.

\begin{theorem}
\label{Pconvex}Let $C$ be a nonempty bounded closed convex subset of a $p$%
-uniformly convex Banach space $X$ with the Opial property and $\mathcal{T}%
=\{T_{t}:t\in G\}$ an asymptotically regular semigroup on $C$ such that%
\begin{equation*}
\liminf_{t}|T_{t}|<\max \left\{ (1+c_{p})^{1/p},\left( \frac{1}{2}\left(
1+(1+4c_{p}WCS(X)^{p})^{1/2}\right) \right) ^{1/p}\right\} .
\end{equation*}%
Then $\mathcal{T}$ has a fixed point in $C$ and $\Fix\mathcal{T}$ is a H\"{o}%
lder continuous retract of $C.$
\end{theorem}

\begin{proof}
Choose a sequence $(t_{n})$ of elements in $G,$ $\lim_{n\rightarrow \infty
}t_{n}=\infty ,$ such that $s(\mathcal{T})=\lim_{n\rightarrow \infty
}\left\vert T_{t_{n}}\right\vert $ and let $(t_{\alpha })_{\alpha \in \emph{A%
}}$ denotes a pointwise weakly convergent subnet of $(t_{n}).$ Define, for
every $x\in C$,%
\begin{equation*}
Lx=w\text{-}\lim_{\alpha }T_{t_{\alpha }}x.
\end{equation*}%
Fix $x_{0}\in C$ and put $x_{m+1}=Lx_{m},m\geq 0.$ Let $B_{m}=\limsup_{%
\alpha }\left\Vert T_{t_{\alpha }}x_{m}-x_{m+1}\right\Vert .$ Since $X$
satisfies the Opial property, it follows from \cite[Prop. 2.9]{KaPr} that
\begin{equation*}
\limsup_{\alpha }\left\Vert T_{t_{\alpha }}x_{m}-x_{m+1}\right\Vert
<\limsup_{\alpha }\left\Vert T_{t_{\alpha }}x_{m}-y\right\Vert
\end{equation*}%
for every $y\neq x_{m+1},$ i.e., $x_{m+1}$ is the unique point in the
asymptotic center $A(C,(T_{t_{\alpha }}x_{m})),m\geq 0.$ Applying (\ref%
{ineq_Xu}) yields%
\begin{align*}
& c_{p}W_{p}(\lambda )\left\Vert x_{m}-T_{t_{\alpha }}x_{m}\right\Vert
^{p}+\left\Vert \lambda x_{m}+(1-\lambda )T_{t_{\alpha }}x_{m}-T_{t_{\beta
}}x_{m-1}\right\Vert ^{p} \\
& \leq \lambda \left\Vert x_{m}-T_{t_{\beta }}x_{m-1}\right\Vert
^{p}+(1-\lambda )\left\Vert T_{t_{\alpha }}x_{m}-T_{t_{\beta
}}x_{m-1}\right\Vert ^{p}.
\end{align*}%
for every $\alpha ,\beta \in \emph{A},0<\lambda <1,m>0.$ Following \cite[Th.
9]{GoTai} (see also \cite{Xu90}) and using the asymptotic regularity of $%
\mathcal{T},$ we obtain%
\begin{equation}
\limsup_{\alpha }\left\Vert T_{t_{\alpha }}x_{m}-x_{m}\right\Vert ^{p}\leq
\frac{s(\mathcal{T})^{p}-1}{c_{p}}(B_{m-1})^{p}.  \label{in1}
\end{equation}%
for any $m>0.$ By Theorem \ref{Wi1} and the weak lower semicontinuity of the
norm, we have
\begin{equation}
B_{m}\leq \frac{D[(T_{t_{\alpha }}x_{m})]}{\WCS(X)}\leq \frac{s(\mathcal{T})%
}{\WCS(X)}\limsup_{\alpha }\left\Vert T_{t_{\alpha }}x_{m}-x_{m}\right\Vert .
\label{in2}
\end{equation}%
Furthermore, by the Opial property,
\begin{equation}
B_{m}\leq \limsup_{\alpha }\left\Vert T_{t_{\alpha }}x_{m}-x_{m}\right\Vert .
\label{in3}
\end{equation}%
Combining (\ref{in1}) with (\ref{in2}) and (\ref{in3}) we see that%
\begin{equation*}
(B_{m})^{p}=\limsup_{\alpha }\left\Vert T_{t_{\alpha
}}x_{m}-x_{m+1}\right\Vert ^{p}\leq \gamma ^{p}(B_{m-1})^{p},
\end{equation*}%
where
\begin{equation*}
\gamma ^{p}=\max \left\{ \frac{s(\mathcal{T})^{p}-1}{c_{p}},\frac{s(\mathcal{%
T})^{p}-1}{c_{p}}\left( \frac{s(\mathcal{T})}{\WCS(X)}\right) ^{p}\right\}
<1,
\end{equation*}%
by assumption. Hence $B_{m}\leq \gamma B_{m-1}$ for every $m\geq 1$ and,
proceeding in the same way as in the proof of Theorem \ref{Thwcs}, we
conclude that $\Fix\mathcal{T}$ is a nonempty H\"{o}lder continuous retract
of $C.$
\end{proof}

\end{document}